\documentclass{agtart_a}
\pdfoutput=1

\title{The homology of the Milnor fiber for classical braid groups}

\author{Filippo Callegaro}
\givenname{Filippo}
\surname{Callegaro}
\address{Scuola Normale Superiore\\Piazza dei Cavalieri, 7\\\newline
56126 Pisa\\Italy}
\email{f.callegaro@sns.it}
\urladdr{}

\volumenumber{6}
\issuenumber{}
\publicationyear{2006}
\papernumber{67}
\startpage{1903}
\endpage{1923}

\doi{}
\MR{}
\Zbl{}


\received{30 November 2005}
\revised{22 September 2006}
\accepted{23 September 2006}
\published{14 November 2006}
\publishedonline{14 November 2006}
\proposed{}
\seconded{}
\corresponding{}
\editor{CPR}
\version{}

\arxivreference{math.AT/0511453}

\AtBeginDocument{\let\bar\wbar\let\tilde\wtilde\let\hat\what}
\def\strutt{\vrule width 0pt depth 6pt height 11pt}

\newcommand{\N}[0]{\mathbb{N}}

\newcommand{\ph}[0]{\varphi}
\newcommand{\bd}[1]{\textbf{#1}}

\newcommand{\p}[0]{p}

\newcommand{\sst}[0]{\subset}
\renewcommand{\cl}[1]{\mathcal{#1}}

\newcommand{\into}[0]{\hookrightarrow}

\newcommand{\im}[0]{\mbox{im }}
\newcommand{\lar}[1]{\stackrel{#1}{\leftarrow}}
\newcommand{\de}[0]{\partial}
\newcommand{\api}[0]{\p^{i-1}}
\newcommand{\bpi}[0]{\phi(p^i)}
\renewcommand{\S}[0]{\mathfrak{S}}

\newcommand{\pmu}[0]{{\pm 1}}

\newcommand{\Br}[0]{\mathrm{Br}}
\newcommand{\A}[0]{\mathrm{A}}

\makeatletter
\def\cnewtheorem#1[#2]#3{\newtheorem{#1}{#3}[section]
\expandafter\let\csname c@#1\endcsname\c@df}
\makeatother

\theoremstyle{definition}

\theoremstyle{plain}
\newtheorem{teo}{Theorem}
\makeautorefname{teo}{Theorem}
\newtheorem*{teo*}{Theorem}
\cnewtheorem{prop}[df]{Proposition}
\cnewtheorem{lem}[df]{Lemma}
\cnewtheorem{cor}[df]{Corollary}

\begin{document}

\begin{asciiabstract}
In this paper we compute the homology of the braid groups, with
coefficients in the module Z[q^+-1] given by the ring of Laurent
polynomials with integer coefficients and where the action of the
braid group is defined by mapping each generator of the standard
presentation to multiplication by -q.

The homology thus computed is isomorphic to the homology with constant
coefficients of the Milnor fiber of the discriminantal singularity.
\end{asciiabstract}

\begin{htmlabstract}
<p class="noindent">
In this paper we compute the homology of the braid groups, with
coefficients in the module <b>Z</b>[q<sup>+-1</sup>] given by the ring of
Laurent polynomials with integer coefficients and where the action of
the braid group is defined by mapping each generator of the standard
presentation to multiplication by -q.
</p>
<p class="noindent">
The homology thus computed is isomorphic to the homology with constant
coefficients of the Milnor fiber of the discriminantal singularity.
</p>
\end{htmlabstract}

\begin{abstract}
In this paper we compute the homology of the braid groups, with
coefficients in the module $\mathbb{Z}[q^{\pm1}]$ given by the ring of
Laurent polynomials with integer coefficients and where the action of
the braid group is defined by mapping each generator of the standard
presentation to multiplication by $-q$.

The homology thus computed is isomorphic to the homology with constant
coefficients of the Milnor fiber of the discriminantal singularity.
\end{abstract}

\maketitle
\section{Introduction}
Let $(W,S)$ be a Coxeter system, with $W$ a finite, irreducible
Coxeter group and let $G_W$ be the associated Artin group (see
Bourbaki \cite{bourb} for an introduction to Coxeter groups and their
classifications and Brieskorn and Saito \cite{brs} for relations
between Coxeter groups and Artin groups ). The main objects of study
of this paper are the Artin groups of type $\A_n$. We recall that the
Artin group $G_{\A_n}$ is the same as the classical braid group (see
Artin \cite{art1} for the original definition and Birman \cite{bir}
for more details) on $n+1$ strands $\Br(n+1)$ and a standard
presentation is given by:
$$
 \langle \sigma_i, i=1, \ldots, n \mid \sigma_i \sigma_{i+1} \sigma_i = \sigma_{i+1} \sigma_i \sigma_{i+1}, \sigma_i \sigma_j = \sigma_j \sigma_i \mbox{ if } | i-j | \geq 2 \rangle.
$$
The homology and cohomology with trivial coefficients for Artin groups
$G_W$ associated to Coxeter groups $W$ of finite type are well
known. The first partial computations for the cohomology of braid
groups are due to Arnol'd and appear in \cite{ar} and \cite{ar2}. In
\cite{fuks} Fuks computes the cohomology ring of braid groups with
$\Z_2$--coefficients using a cell decomposition of the Alexandroff
compactification of the configuration spaces of $\C$, which are
classifying spaces for braid groups. Using the same cell
decomposition, in \cite{vain} Va{\u\i}n{\v{s}}te{\u\i}n computes the
cohomology ring of braids with $\Z_p$--coefficients for any prime $p$;
moreover he computes the Bockstein operator, hence he gets the
cohomology ring of braids with $\Z$--coefficients. The same results in
homology are obtained independently by F Cohen in \cite{coh}: he
computes the homology of braid groups over $\Z$ using the theory of
homology operations in $n$--fold loop spaces. In \cite{gorj} Gorjunov
extends the results of Va{\u\i}n{\v{s}}te{\u\i}n computing the
cohomology ring of Artin groups of type ${\mathrm C}_n$ and ${\mathrm
D}_n$ over $\Z$. Finally in \cite{salv2} Salvetti computes the
cohomology groups with $\Z$--coefficients for the Artin groups
associated with exceptional Coxeter groups; the ring structure is
determined by Landi in \cite{lan}.

Let $\bd{X}_W$ be the classifying space for $G_W$. For a ring $A$ we
can define the local system $A[q^\pmu]$ over the space $\bd{X}_W$,
with twisted coefficients over the ring $A[q^\pmu]$, where each
standard generator of $G_W$ acts as multiplication by $-q$.

The homology groups $H_*(\bd{X}_W; A[q^\pmu])$ are equal to the
homology group (with trivial coefficients over the ring $A$) of the
Milnor fiber $\bd{F}_W$ of the discriminant singularity associated to
$W$.  The same hold for the cohomology groups, modulo an index
shifting ((as proved by the author in \cite{c1}). Hence we get a
$q$--module structure on the homology and cohomology of the Milnor
fiber $\bd{F}_W$, where the $-q$--multiplication corresponds to the
map induced by the monodromy automorphism of the fiber. We remark also
that for the case $\A_n$ the fiber $\bd{F}_{\A_n}$ is a classifying
space for the commutator subgroup of the braid group $\Br(n+1)$, hence
the homology (cohomology) of the fiber is also the homology
(cohomology) of $\Br(n+1)'$.

The groups $H^*(\bd{X} _W; A[q^\pmu])$ have been computed by the
author and Salvetti \cite{cal-sal} for all exceptional Artin
groups. Several authors (De Concini, Procesi and Salvetti \cite{dps},
Frenkel' \cite{fren}, Markaryan \cite{mar}) made this computation for
type $\A_n$ Artin groups, with $A$ a field of characteristic $0$. In
de Concini, Procesi, Salvetti and Stumbo \cite{dpss} the computations
have been performed for all other finite type Artin groups, with
rational coefficients. In all these cases all cohomology modules are
sums of cyclotomic fields (or zero). In particular, in \cite{dps}, an
interesting arithmetic behaviour of the table $(H^i(\Br(j);
A[q^\pmu]))_{i,j}$ is noticed.

In \cite{mar} Markaryan used the isomorphism between the standard
resolution of a certain algebra and the algebraic complex associated
to the classifying spaces for braid groups to compute the homology of
braid groups with coefficients in the local system ${\Q[q^\pmu]}$.

In this paper we extend the use of this resolution in order to compute
the homology of braid groups with coefficients in the local system
${K[q^\pmu]}$ for a generic field $K$. We also compute the Bockstein
operator in order to get the homology over ${\Z[q^\pmu]}$. Some
computations for the first cases can be found in tables \ref{ta:1},
\ref{ta:2}, \ref{ta:3} and \ref{ta:4}. Our main result is stated in
\fullref{t:summary}.

In \cite{cp} Cohen and Pakianathan compute the homology of the braid group on infinitely many strands $\Br(\infty)$ with coefficients in the local system $K[q^\pmu]$ for any field $K$: this is the stable part of the homology of $\Br(n)$ (with coefficients the same local system) with respect to the embeddings $j_n\co \Br(n) \into \Br(n+1)$.
In \fullref{c:stable} we obtain the same result; moreover we are able to compute the Bockstein operator, hence we give a presentation of the homology of $\Br(\infty)$ with coefficients in the local system $\Z[q^\pmu]$.

\medskip
\textbf{\fullref{c:stable}}\qua
\sl We have that
$
H_*(\Br(\infty); \Q[q^\pmu])= \Q,
$
concentrated in dimension $0$;
$$
H_*(\Br(\infty); \Z_2[q^\pmu])= \Z_2 [x_2^2, x_{2^i}, i>1]
$$
and for a prime $p > 2$
$$
H_*(\Br(\infty); \Z_p[q^\pmu])= \Z_p 
[y_{2p^i}, x_{2p^i}, i > 0 ]/(x_{2p^{i}}^2)
$$
with $\dim x_i = i-1$, $\dim y_i = i-2$. The Bockstein operator acts as follows:
$$
\beta_2 x_{2^i} = x_{2^{i-1}}^2;\qquad\beta_p y_i = 0; 
\qquad \beta_p x_i = y_i \qquad \mbox{(for }p>2 \mbox{)}.
$$
The homology $H_*(\Br(\infty); \Z[q^\pmu])$ has no $p^2$--torsion for any prime $p$. A presentation of $H_*(\Br(\infty); \Z[q^\pmu])$ is given by
$$
\Z\left[
\begin{array}{c}
y_{2p^i},x_{2^j}^2,\\
x_{2^i}^2x_{2^{i_1}} \cdots x_{2^{i_h}} \\
y_{2p^j} x_{2p^{j_1}} \cdots x_{2p^{j_h}}
\end{array} 
\right] \left/ \left( 
\begin{array}{c}
2 x_{2^i}, py_{2p^j}, x_{2p^j}^2
\end{array}
\right) \right.
$$
with indices running as follows: $0<i$,$i+1 < i_1 < \cdots < i_h$, $0<j< j_1 < \cdots < j_h$ and $p$ in the set of odd primes. The structure of $\Z[q^\pmu]$--module is trivial and so the action of $q$ corresponds to multiplication by $-1$. 
\rm \medskip

\section{Notation and definitions} \label{s:due}

Let $R$ be a ring with identity and let $q$ be an element of $R$.
Following \cite{mar} we define the \emph{algebra of $q$--divided polynomials} $\Gamma_{R}(t,q)$ as the graded algebra over $R$ with  
generators $t_i$ ($ i \in \N$, $\deg t_i = i$) and relations
\begin{gather*}
t_it_j = \left[ 
\begin{array}{c}
i+j \\ i 
\end{array}
\right]t_{i+j}\\
\left[ 
\begin{array}{c}
i+j \\ i 
\end{array}
\right]\vrule width 0pt depth 15pt = \frac{(1-q)(1-q^2)\cdots(1-q^{i+j})}{(1-q)\cdots(1-q^i)(1-q)\cdots(1-q^j)}\tag*{\hbox{where}}
\end{gather*}
is the \emph{$q$--binomial} coefficient; if we define the \emph{$q$--analog} of $i$ as  
$$ [i] = 1+q+\cdots +q^{i-1}, $$ 
and $[i]!= [2][3]\cdots[i]$, then we can write also 
$$
\left[ 
\begin{array}{c}
i+j \\ i 
\end{array}
\right] = \frac{[i+j]!}{[i]![j]!}. 
$$
We recall that if $q$ commutes with $a$ and $b$ and $ba=qab$, then
$$
(a+b)^n = \sum_{i=0}^n\left[ 
\begin{array}{c}
n \\ i 
\end{array}
\right] a^ib^{n-i}.
$$
Now we want to study the homology and cohomology (as defined in Cartan
and Eilenberg \cite{ce}) of the algebra $\Gamma_{R}(t,q)$. We can
consider the \emph{normalized standard complex} (see \cite{ce} for a
general definition) that calculates the homology of the algebra
$\Gamma_{R}(t,q)$. The complex is given as follows:
$$
0 \lar{} R = C_0 \lar{\de} C_1 \lar{\de} C_2 \lar{\de} \cdots ,
$$
where the $R$--module $C_n$  is freely generated by all the monomials of the form $a\otimes t_{i_1} \otimes \cdots \otimes t_{i_n}$, with $a \in R$ and the boundary formula is:
\begin{eqnarray*}
&   & \de(a\otimes t_{i_1} \otimes \cdots \otimes t_{i_n}) = \\
& = & \sum_{j=1}^{n-1} (-1)^{j+1} a\otimes t_{i_1} \otimes \cdots \otimes t_{i_j} t_{i_{j+1}} \otimes \cdots \otimes t_{i_n}  = \\
& = & \sum_{j=1}^{n-1} (-1)^{j+1} \left[ \begin{array}{c} i_j+i_{j+1} \\ i_j \end{array} \right] a\otimes t_{i_1} \otimes \cdots \otimes t_{i_j+i_{j+1}} \otimes \cdots \otimes t_{i_n}. 
\end{eqnarray*}
By means of the grading of the algebra $\Gamma_{R}(t,q)$, the module $C$ is decomposed into the direct sum of complexes of different degrees,
$$
C = \bigoplus_{i=0}^\infty C^{(i)}, 
$$
where $\deg(a\otimes t_{i_1} \otimes \cdots \otimes t_{i_n})= i_1 + \cdots + i_n$ and for $c \in C_k^{(n)}$ we set $\deg c = n$ and $\dim c = k$.

The dual complex $C^*$, given by the modules $C^n= \mbox{Hom}(C_n, R)$ and with coboundary map the transposed map of $\de$, computes the cohomology ring of the algebra $\Gamma_{R}(t,q)$. The multiplication is defined on representatives as follows: if $m^*_1$ and $m^*_2$ are the dual classes of the monomials $m_1$ and $m_2$, respectively, then the product $m^*_1 m^*_2$ is the dual class of the monomial $m_1 \otimes m_2$.

Given a space $X$ such that $\pi_1(X)= \Br(n)$, we can define a local system $R$ on $X$. Over a point $x \in X$ we have the ring $R$; the system of coefficients is twisted and the action is given by sending each standard generator of the group $\Br(n)$ to multiplication by $-q$. This action corresponds to the determinant of the Burau representation for the braid group $\Br(n)$ (see, for example, \cite{cp}). We remark that although the choice of the multiplication by $q$ would be equivalent, we use $-q$, which  seems more natural to us, and also for coherence with \cite{cp}, \cite{dps}, \cite{dpss}.

The complex $C^{(n)}_{n-*}$ coincides with the complex (defined in \cite{salv2}) that computes the cohomology of the group $\Br(n)$ with local coefficients $R$.

By the module $H_*(\Br(*),R)$ we mean the bigraded module (the gradings are the degree $\deg$ and the dimension $\dim$) whose component of degree $n$ and dimension $l$ is $H_l(\Br(n),R)$. We can think of $H_*(\Br(*),R)$ as a ring using the multiplication induced by the standard homomorphism (obtained by juxtaposing braids)
$$
\mu_{ij}\co \Br(i) \times \Br(j) \to \Br(i+j).
$$

\begin{teo} {\rm\cite{mar}}\qua
The ring $H_*(\Br(*),R)$ coincides, modulo a change of indexes, with the cohomology ring of the algebra $\Gamma_R(t,q)$:
$$
H_l(\Br(n), R) \simeq H^{n-l}{(\Gamma_R(t,q))}_{(\deg = n)}
$$ 
and the product structure in $H_*(\Br(*),R)$ coincides with the cohomological multiplication of the ring $H^*(\Gamma_R(t,q))$. 
\end{teo}

\section{The Milnor fiber and some lemmas} \label{s:milnor}
Let $V =\C^n$ be a finite-dimensional complex vector space. The symmetric group on $n$ elements  $\S_n$ acts on this space by permuting the coordinates. Let $l_{ij}$ be the linear functional $z_i - z_j$ and let $\cl{H}_{ij}$ be the hyperplane $\{ l_{ij} = 0 \} $. The complement of the union of the hyperplanes
$$
\bd{Y}_n = V \setminus \bigcup_{i<j}\cl{H}_{ij}
$$
is a classifying space for the pure braid group on $n$ strands. If we consider the quotient of $\bd{Y}_n$ with respect to the action of $\S_n$
$$
\bd{X}_n = \bd{Y}_n/\S_n
$$
we get a classifying space for the braid group $\Br(n)$. Consider the product $\delta = \prod_{i<j} l_{ij}^2$. The polynomial $\delta$ is invariant with respect to the action of $\S_n$ and so it induces a map 
$$
\delta'\co \bd{X}_n \to \C^*.
$$
The fiber $\bd{F}_1(n) = {\delta'}^{-1}(1)$ is the Milnor fiber of the discriminant singularity $\bd{F}_0(n) = \bigcup \cl{H}_{ij} / \S_n$ in the affine variety $V/\S_n$ (which is isomorphic to the complex space $\C^n$). The complement of $\bd{F}_0(n)$ in $V/\S_n$ can also be thought as the set of polynomials with distinct roots in the space of all monic polynomials of degree $n$ with complex coefficients. 
Moreover, the fiber $\bd{F}_1(n)$ is a classifying space for the commutator subgroup $\Br(n)'$ of the braid group $\Br(n)$ and we have that (see for example \cite{c1}):
\begin{gather*}
H_*(\Br(n)',\Z) \simeq H_*(\bd{F}_1(n),\Z) \simeq H_*(\Br(n),\Z[q^\pmu]) 
\\
H^*(\Br(n)',\Z) \simeq H^*(\bd{F}_1(n),\Z) \simeq H^{*+1}(\Br(n),\Z[q^\pmu]),
\tag*{\hbox{and}}
\end{gather*}
with the usual $q$--action.

In what follows $K$ is a field and $p$ always refers to the  characteristic of the field $K$ ($p=0$ or $p$ a prime).  
Cyclotomic polynomials are usually defined over a field of characteristic $0$, by saying that the $m$-th cyclotomic polynomial is the monic polynomial whose roots are all simple roots and are all the $m$-th primitive roots of unity. Over a generic field $K$ we can define by induction the $m$-th cyclotomic polynomial $\ph_m$, by saying that $\ph_1=q-1$ and $q^m-1 = \prod_{i \mid m} \ph_i$. For each positive integer $m$ we define the ring $K(m) =  K[q]/\ph_m$.
\vspace{3pt}

We have the following technical lemmas: 
\vspace{3pt}

\begin{lem} {\rm\cite{c}}\qua\label{l:tecnico1}
Let $m < n$ be two positive integers.
Then we have:
$$
(\ph_m, \ph_n) = \left\{ \begin{array}{cl}
(\ph_m, p) & \mbox{if }n=mp^i, i \geq1, \mbox{for a prime } p \\
(1) & \mbox{otherwise.}
                         \end{array} \right.
$$ 
\end{lem}

\vspace{3pt}
We leave the proof to the Appendix.
As an easy consequence of this Lemma we obtain the following Corollary, whose proof is left to the reader:
\vspace{3pt}
\begin{cor} \label{c:tecnico}
Let $i<j$ be two positive integers. Then we can write
\begin{equation} \label{eq:tecn1}
\ph_{p^j} = \ph_{p^i} \omega + p \psi
\end{equation} 
where $\omega, \psi \in \Z[q^\pmu]$ and $\psi$ is invertible $\mod \ph_{p^i}$;
\begin{equation} \label{eq:tecn2}
\ph_{mp^j} = \ph_{mp^i} \omega + p \psi
\end{equation} 
where $\omega, \psi \in \Z[q^\pmu]$ and $\psi$ is invertible $\mod \ph_{mp^i}$. 
\end{cor} 
\vspace{3pt}
We can fix once and for all polynomials $\omega_{p^j, p^i}$, $\omega_{mp^j, mp^i}$, $\psi_{p^j, p^i}$, $\psi_{mp^j, mp^i}$ satisfying the equations \eqref{eq:tecn1} and \eqref{eq:tecn2}.
\vspace{4pt}

\begin{lem}\label{l:tecnico2}{\rm\cite{guer}}\qua 
Let $m$ be an integer and $p$ a prime. Then we have: 
\begin{equation} \label{prima}
\ph_{p^i} \simeq \ph_p^{p^{i-1}}\quad  \mod p;
\end{equation} 
if we suppose that $p\nmid m$, then: 
\begin{equation} %\label{seconda}
\ph_{mp^i} \simeq \ph_m^{\phi(p^i)}\quad  \mod p \label{prima_equiv}
\end{equation}
where $\phi$ denotes the Euler $\phi$--function. 
\end{lem}

Now we consider again the algebra of $q$--divided polynomials $\Gamma_R(t,q)$ in the case $R = K(m)$. 

\begin{lem}
The following decompositions hold:

{\rm (a)}\qua if $p=0$ (see also \cite{mar}):
\begin{equation} 
\Gamma_{K(m)}(t,q) \simeq K(m)[u_m] \otimes K(m)[u_1]/(u_1^m);
\end{equation} 

{\rm (b)}\qua if $p \neq 0$:
\begin{equation}
\Gamma_{K(p)}(t,q) \simeq \bigotimes_{i=0}^\infty K(p)[u_{p^i}]/(u_{p^i}^p); 
\end{equation} 

{\rm (c)}\qua if $p \neq 0$ and $p \nmid m$:
\begin{equation}
\Gamma_{K(m)}(t,q) \simeq K(m)[u_1]/(u_1^m) \otimes \bigotimes_{i=0}^\infty K(m)[u_{p^im}]/(u_{p^im}^p); 
\end{equation} 
with $\deg u_j =j$.
\end{lem}

\begin{proof}
The proof of (a) is given in \cite{mar}. 

For (b) the isomorphism is given as follows:
\begin{equation} \label{eq:pexpress}
t_i \mapsto \frac{u_1^{k_0}u_p^{k_1} \cdots u_{p^r}^{k_r}}{[i]!/\ph_p^{{[i]!}_{\ph_p}} }
\end{equation}
where $k_r \cdots k_0$ is the expression of $i$ in base $p$ and ${x}_{\ph_p}$ is the maximal exponent $l$ such that $\ph_p^l$ divides $x$.

For (c) we have the isomorphism given by
\begin{equation} \label{eq:mpexpress}
t_i \mapsto \frac{u_1^ku_m^{k_0}u_{mp}^{k_1} \cdots u_{mp^r}^{k_r}}{[i]!/\ph_m^{{[i]!}_{\ph_m}}}
\end{equation}
where $k$ is the remainder of the division of $i$ by $m$ and $k_r \cdots k_0$ is the expression of $(i-k)/m$ in base $p$.
 
The Lemma follows from the next key observation: if $k_r \cdots k_0$ is the expression of $i$ in the base $p$ and $k_r' \cdots k_0'$ is the expression for $j$ (resp. $k$, $k_r \cdots k_0$ and $k'$, $k_r' \cdots k_0'$ are the numbers associated to $i$ and $j$ as in \eqref{eq:mpexpress}), then the polynomial $\ph_p$ (resp. $\ph_m$) does not divide 
$ \left[ \begin{array}{c} i+j \\ i \end{array} \right] $ if and only if the expression for $i+j$ in base $p$ is given by $h_r \cdots h_0$, with $h_l = k_l + k_l'$ for $l=0, \ldots, r$ (resp. the numbers associated to $i+j$ are $h$, $h_r \cdots h_0$, with $h =k+k'$, $h_l = k_l + k_l'$ for $l=0, \ldots, r$).
\end{proof}

The cohomology rings of $R[u]$ and $R[u]/(u^i)$ are already known. In fact we have:

\begin{lem}{\rm\cite{mar}}\label{l:marl}\quad
$
H^*(R[u]) \simeq 
\begin{array}{cr}
\Lambda[x], & \deg(x)= \deg(u), \dim(x) = 1;
\end{array}
$
$$
H^*(R[u]/(u^n)) \simeq \left\{ 
\begin{array}{cr} 
R[x]  &\mbox{for } n=2\\ 
R[y] \otimes \Lambda[x] & \mbox{for } n>2
\end{array}   \right.
$$
where $\deg(x) = \deg(u) $, $\deg(y) = n \deg(u) $, $\dim(x) = 1 $, $\dim(y) = 2  $ and $\Lambda[x]$ is the exterior algebra over the ring $R$ in the variable $x$. 
\end{lem}

We remark that generators of the rings in the Lemma can be given as follows: a representative $x$ is given by the dual class of $u$. 
Moreover in characteristic $p=0$, a representative of $y$ is given by
$$
\sum_{i=1}^{n-1}\left( \begin{array}{c}
n \\ i \end{array}\right) (u^i \otimes u^{n-i})^*
$$ 
and with $p \neq 0$, if $n$ is a power of $p$, we can choose as a representative
$$
\frac{1}{p} \sum_{i=1}^{n-1}\left( \begin{array}{c}
n \\ i \end{array}\right) (u^i \otimes u^{n-i})^*,
$$ 
where the notation $(u^i \otimes u^{n-i})^*$ means the dual class of $(u^i \otimes u^{n-i})$.

\section{Computations and results}

Now we can calculate the cohomology of $\Gamma_{K_p}(t,q)$ (and so the homology of $\Br(*)$ with coefficients in the local system ${K(m)}$) applying the fact that the cohomology of a tensor product of algebras is the tensor product of the cohomology of the factors. 

Applying \fullref{l:marl} we have the following straightforward results:

\begin{teo}{\rm\cite{mar}}\label{t:mardecompos}\qua
If $p=0$ and $m = 2$ then
\begin{gather*}
H_*(\Br(*); K(2)) \simeq \Lambda[x_2] \otimes K(2)[x_1];
\\
H_*(\Br(*); K(m)) \simeq \Lambda[x_m] \otimes K(m)[y_m] \otimes \Lambda[x_1];
\tag*{\hbox{for $m > 2$:}}
\end{gather*}
with $\deg x_i = i$, $\dim x_i = i-1$, $\deg y_m=m$, $\dim y_m =m-2$. 
\end{teo}
\vspace{4pt}

\begin{teo} \label{t:decompos}
Let $p$ be a prime and $m$ be a positive integer, such that $p \nmid m$. 
We have the following cases:

{\rm(a)}\qua if $p = 2$:\qquad
$
H_*(\Br(*); K(2)) \simeq \bigotimes_{i=0}^\infty K(2)[x_{2^i}];
$
$$
H_*(\Br(*); K(m)) \simeq K(m)[y_m] \otimes \Lambda[x_1] \otimes \bigotimes_{i=0}^\infty K(m)[x_{m2^{i}}];
$$

{\rm(b)}\qua if $p > 2$ and $m = 2$:
\begin{gather*}
H_*(\Br(*); K(p)) \simeq \bigotimes_{i=0}^\infty (K(p)[y_{p^{i+1}}] \otimes \Lambda[x_{p^i}]);
\\
H_*(\Br(*); K(2)) \simeq K(2)[x_1] \otimes \bigotimes_{i=0}^\infty (K(2)[y_{2p^{i+1}}] \otimes \Lambda[x_{2p^i}]);
\end{gather*}

{\rm(c)}\qua if $p>2$ and $m > 2$:
$$
H_*(\Br(*); K(m)) \simeq K(m)[y_m] \otimes \Lambda[x_1] \bigotimes_{i=0}^\infty (K(m)[y_{p^{i+1}m}] \otimes \Lambda[x_{p^im}]);
$$
where $\deg x_i = \deg y_i = i$, $\dim x_i = i- 1$, $\dim y_i = i - 2$. 
\end{teo}
\vspace{4pt}

We want to use these results to compute the homology of $\Br(*)$ with coefficients in the local system $A$ over the ring $A = K[q^\pmu]$ with the same twisting defined as in \fullref{s:due}. 
\vspace{4pt}

The exact sequence
$$
1 \to \Br(n)' \into \Br(n) \to \Z \to 1,
$$
tells us that the homology $H_*(\Br(n); A)$ is $H_*(\Br(n)'; K)$ as an $A$--module (see for example \cite{mar}, \cite{coh-suc} or \cite{c1}); since for $n \neq 3,4$, $\Br(n)'' = \Br(n)'$ (see \cite{gl} for a proof of this), we have that $H_0(\Br(n)'; K) = K$, $H_1(\Br(n)';K)=0$. Moreover the $A$--action on $H_0$ is trivial and so $H_0(\Br(n); A)= A/(q+1)$ as an $A$--module. Moreover we have:
\vspace{4pt}

\begin{lem}{\rm\cite{mar}}\label{l:torsion}\qua
The $R$--modules $H_l(\Br(n),A) (n > 1, l> 0)$ are annihilated by multiplication by $[n]!$. 
\end{lem}
\vspace{3pt}

Let us consider a polynomial $a \in A$. We can consider the set $S_a$ of all elements $b \in A$ that are prime with $a$. It is clear that $S_a$ is a multiplicatively closed set. We write  
$A_{(a)}$ for the localization $A_{S_a}$ of the ring $A$ respect to the set $S_a$. 
\vspace{3pt}

It follows from \fullref{l:tecnico1} that for $p =0$, $\ph_m$ is invertible in $A_{(\ph_n)}$ if and only if $m \neq n$; for $p \neq 0$, $\ph_m$ is invertible in $A_{(\ph_n)}$ if and only if $n \neq mp^i$ and $m \neq np^i$, $\forall i \geq1$.
\vspace{3pt}

The following decompositions hold for the homology of $\Br(n)$ with coefficients in the local system ${A}$: 
\vspace{3pt}

\begin{lem} \label{l:decompos}
Let $n > 1$. For $p = 0$ we have: 
\begin{gather*}
H_*(\Br(n); A) \simeq 
\bigoplus_{ m = 2}^\infty
H_* (\Br(n); {A_{(\ph_m)}});
\\
H_*(\Br(n); A) \simeq 
\bigoplus_{ p \nmid m \mbox{\begin{scriptsize} or \end{scriptsize}} m=p} 
H_* (\Br(n); {A_{(\ph_m)}}).
\tag*{\hbox{for $p \neq 0$:}}
\end{gather*}
\end{lem}
\vspace{3pt}

\begin{proof}
Consider the homomorphism 
$$
{i_m}_*\co H_*(\Br(n); A) \to H_* (\Br(n); {A_{(\ph_m)}})
$$
induced by the injection $i_m\co A \into A_{(\ph_m)}$. 
We extend in a natural way the map ${i_m}_*$ through the tensor product with $A_{(\ph_m)}$ and we get the new map
\begin{equation} \label{e:itilda}
\widetilde{i_m}\co H_*(\Br(n); A) \otimes_A A_{(\ph_m)} \to 
H_*(\Br(n);{A_{(\ph_m)}}).
\end{equation} 
Using Lemmas \ref{l:tecnico1}, \ref{l:tecnico2} and \ref{l:torsion} it is easy to see that in order to prove \fullref{l:decompos} it is enough to show that the map $\widetilde{i_m}$ is an isomorphism.  

First we prove the injectivity of $\widetilde{i_m}$. 
Let $\alpha$ be a representative of an element $v$ in $H_*(\Br(n); A) \otimes_A A_{(\ph_m)}$.
If the corresponding class of $\widetilde{i_m}v$ is zero, and so $\widetilde{i_m} \alpha$ is a boundary, then there exists an element $\beta$ such that $d\beta = \widetilde{i_m} \alpha$. Multiplying $\beta$ by an appropriate polynomial $\psi$ prime with $\ph_m$, we get an element $\beta' = \psi \beta$ that belongs to the resolution of $\Br(n)$ over $A$, so $d\beta' = \psi \alpha$. This means that $\psi \alpha$ belongs to the zero class in $H_*(\Br(n); A)$ and, since $\psi$ is invertible in $A_{(\ph_m)}$, $\alpha$ belongs to the zero class in $H_*(\Br(n); A) \otimes_A A_{(\ph_m)}$. This proves the injectivity of $\widetilde{i_m}$. 

To prove the surjectivity of $\widetilde{i_m}$ we consider a class $w$ in $H_* (\Br(n); {A_{(\ph_m)}})$ and we choose a representative $\beta$ for $w$. Multiplying $\beta$ by an appropriate polynomial $\theta$ prime with $\ph_m$ we get an element $\beta' = \theta \beta$ in the resolution for $H_*(\Br(n); A)$ and we have that 
$$
\widetilde{i_m}(\beta' \otimes \theta^{-1}) = \beta.
$$
This completes the proof. 
\end{proof}

The next step is to compute $H_*(\Br(n);{A_{(\ph_m)}})$. To do this, consider the following short exact sequence:
$$
0 \to A_{(\ph_m)} \stackrel{\ph_m}{\into} A_{(\ph_m)} \stackrel{\pi}{\to} K(m) \to 0
$$
where the first map is multiplication by $\ph_m$.
We want to study the corresponding homology long exact sequence: 
\begin{eqnarray*}
& \cdots & \stackrel{\pi_*}{\to} H_{l+1}(\Br(*);K(m)) \stackrel{\beta}{\to} \\
\stackrel{\beta}{\to} H_l(\Br(*);A_{(\ph_m)}) \stackrel{{(\ph_m)}_*}{\to} & H_l(\Br(*);A_{(\ph_m)}) & \stackrel{\pi_*}{\to} H_l(\Br(*);K(m)) \stackrel{\beta}{\to} \\
\stackrel{\beta}{\to} H_{l-1}(\Br(*);A_{(\ph_m)}) \stackrel{{(\ph_m)}_*}{\to} & \cdots &
\end{eqnarray*} 
We can decompose $H_l(\Br(*);{A_{(\ph_m)}})$ as a direct sum of terms $A/(\psi^i)$, where $\psi$ is a prime factor of $\ph_m$. So, if $H_l(\Br(*);A_{(\ph_m)})$ has a direct summand $A/(\psi^i)$, generated by an element $v$, it follows that $H_{l+1}(\Br(*);K(m))$ and $H_{l}(\Br(*);K(m))$ have as direct summand a copy of $A/(\psi)$ generated respectively by $w$ and $w'$ and we have that
\begin{gather*}
\beta w = \psi^{i-1} v
\\
\pi_* v = w'.\tag*{\hbox{and}}
\end{gather*}
In \fullref{t:mardecompos} (case $p=0$) we have these maps (see also \cite{mar}):
$$
\beta x_m = \widetilde{y_m}, \pi_* \widetilde{y_m} = y_m, \beta x_1 =0. 
$$
while in \fullref{t:decompos} (case $p \neq 0$), the homomorphisms act as follows:
$$
\beta y_i = 0,\qquad  \beta x_1 = 0, 
$$
$$
\beta x_{2^i}= \ph_2^{2^{i-1}-1} \widetilde{x_{2^{i-1}}}^2; \qquad \beta x_{p^i} = \ph_p^{\api-1}\widetilde{y_{p^i}} \quad \mbox{for }p > 2;
$$
$$
\beta x_{mp^i} = \ph_m^{\phi(p^i)-1} \widetilde{y_{mp^i}}\quad \mbox{for }i>0 \mbox { or } m > 2;
$$
where $\widetilde{x_i}$ and $\widetilde{y_i}$ are defined such that:
$$
\pi_* \widetilde{x_i} = x_i; \pi_* \widetilde{y_i} = y_i.
$$
We can now state the following result:

\begin{prop} \label{p:phitorsione}
For $p=0$ we have, for $m=2$:
\begin{gather*}
H_*(\Br(*);{A_{(\ph_2)}}) \simeq A_{(\ph_2)}[x_1]/(\ph_2 x_1^2);
\\
\tag*{\hbox{and for $m >2$:}}
H_*(\Br(*);{A_{(\ph_m)}}) \simeq A_{(\ph_m)}[x_1, y_m]/(x_1^2, \ph_m y_m);
\end{gather*}
if $ p \neq 0$ and $p \nmid m$ we have the following cases:
\vspace{4pt}

{\rm(a)}\qua for $p=2$: 
$$
H_*(\Br(*);{A_{(\ph_2)}})  \simeq A_{(\ph_2)} \left[
\begin{array}{c} x_1, x_{2^j}^2, \\ x_{2^i}^2 x_{2^{i_1}} \cdots x_{2^{i_h}}
\end{array} 
\right] \left/ \left( 
\begin{array} {c} \ph_2^{2^i} x_{2^i}^2 \end{array}
\right) \right. 
$$
with $0\leq i$, $i+1< i_1< \cdots < i_h$, $0<j$;
$$
H_*(\Br(*);{A_{(\ph_m)}})  \simeq A_{(\ph_m)} \left[
\begin{array}{c} x_1, y_m, x_{m2^i}^2, \\ x_{m2^i}^2 x_{m2^{i_1}} \cdots x_{m2^{i_h}},  \\ 
y_m x_{m2^{j_1}} \cdots x_{m2^{j_h}}
\end{array} 
\right] \left/ \left( 
\begin{array} {c} x_1^2, \ph_m y_m, \\ \ph_m^{2^i}x_{m2^i}^2 \end{array}
\right) \right. 
$$
with $0 \leq i$, $i+1 < i_1< \cdots < i_h$, $0 < j_1 < \cdots < j_h$;
\vspace{4pt}

{\rm(b)}\qua for $p > 2$ and $m = 2$: 
$$
H_*(\Br(*);{A_{(\ph_p)}})  \simeq A_{(\ph_p)} \left[
\begin{array}{c} x_1, y_{p^i}, \\ y_{p^i}x_{p^{i_1}} \cdots x_{p^{i_h}}
\end{array} 
\right] \left/ \left( 
\begin{array} {c} x_{p^i}^2, \\ \ph_p^{\api} y_{p^i} \end{array}
\right) \right. 
$$
with $0 < i< i_1< \cdots < i_h$;
$$ 
H_*(\Br(*);{A_{(\ph_2)}})  \simeq A_{(\ph_2)} \left[
\begin{array}{c} x_1, y_{2p^i}, \\
x_1^2x_{2p^{j_1}}\cdots x_{2p^{j_h}} \\ 
y_{2p^i}x_{2p^{i_1}} \cdots x_{2p^{i_h}}
\end{array} 
\right] \left/ \left( 
\begin{array} {c} \ph_2 x_1^2, x_{2p^i}^2, \\ \ph_2^{\phi(p^i)} y_{2p^i} \end{array}
\right) \right. 
$$
with $0 < i< i_1< \cdots < i_h$, $0<j_1 < \cdots <j_h$;
\vspace{4pt}

{\rm(c)}\qua for $p >2$ and $m > 2$:
$$ 
H_*(\Br(*);{A_{(\ph_m)}})  \simeq A_{(\ph_m)} \left[
\begin{array}{c} x_1, y_{mp^i}, \\ y_{mp^i}x_{mp^{i_1}} \cdots x_{mp^{i_h}}
\end{array} 
\right] \left/ \left( 
\begin{array} {c} x_1^2, x_{mp^i}^2, \\ \ph_m^{\phi(p^i)} y_{mp^i} \end{array}
\right) \right. ;
$$
with $0 \leq i< i_1< \cdots < i_h $;

in all cases  $\deg x_i= \deg y_i = i, \dim x_i = i-1$, $\dim y_i = i-2$.

\end{prop}
\vspace{4pt}

In order to get $H_*(\Br(n);\Z[q^\pmu])$ we still have to compute the Bockstein homomorphism $\beta_p$ associated to the short exact sequence 
$$
0 \to \Z_p \into \Z_{p^2} \to \Z_p \to 0.
$$
We'll see that, as in the case of trivial coefficients (see \cite{coh} or \cite{vain}), there is no $p^2$--torsion in the homology of braid groups. 

In case (a) of \fullref{p:phitorsione} the Bockstein acts as follows (the coefficients $\psi_{i,j}$ are those defined in \fullref{c:tecnico}):
$$
\beta_2 x_1 = 0, \quad \beta_2 x_{2^i}^2 = 0, 
$$
$$
\beta_2 x_{2^i}^2 x_{2^{i_1}} \cdots x_{2^{i_h}} = 
\sum_{s=1}^h \left( \psi_{2^{i_s},2^{i}} x_{2^i}^2 x_{2^{i_s-1}}^2 \prod_{t\neq s} x_{2^{i_t}} \right) ;
$$
$$
\beta_2 y_m= 0, \beta_2 x_{m2^i}^2= 0,
$$
$$ 
\beta_2 x_{m2^i}^2 x_{m2^{i_1}} \cdots x_{m2^{i_h}} =
\sum_{s=1}^h \left( \psi_{m2^{i_s-1},m2^i} x_{m2^i}^2 x_{m2^{i_s-1}}^2 \prod_{t\neq s} x_{m2^{i_t}} \right),
$$
$$
\beta_2 y_m x_{m2^{i_1}} \cdots x_{m2^{i_h}}= 
\sum_{s=1}^h \left( \psi_{m2^{i_s-1},m} y_m x_{m2^{i_s-1}}^2 \prod_{s \neq t} x_{m2^{i_t}} \right);
$$
in case (b) we have:
$$
\beta_p x_1 =0, \beta_p y_{p^i} =0,
$$
$$
\beta_p y_{p^i}x_{p^{i_1}} \cdots x_{p^{i_h}} = 
-\sum_{s=1}^h \left( \psi_{p^{i_s},p^i} y_{p^i}y_{p^{i_s-1}} \prod_{s \neq t} x_{p^{i_t}} \right);
$$
$$
\beta_p x_1^2x_{2p^{i_1}} \cdots x_{2p^{i_h}} =
-\sum_{s=1}^h \left( \psi_{2p^{i_s},2p^i} x_1^2 y_{2p^{i_s-1}} \prod_{s \neq t} x_{2p^{i_t}} \right);
$$
$$
\beta_p y_{2p^i} = 0,
$$
$$
\beta_p y_{2p^i}x_{2p^{i_1}} \cdots x_{2p^{i_h}} = 
-\sum_{s=1}^h \left( \psi_{2p^{i_s},2p^i} y_{2p^i}y_{2p^{i_s-1}} \prod_{s \neq t} x_{2p^{i_t}} \right);
$$
finally in case (c) the map is:
$$
\beta_p x_1 = 0, \beta_p y_{mp^i} = 0,
$$
$$
\beta_p y_{mp^i}x_{mp^{i_1}} \cdots x_{mp^{i_h}} = 
-\sum_{s=1}^h \left( \psi_{mp^{i_s},mp^i} y_{mp^i}y_{mp^{i_s-1}} \prod_{s \neq t} x_{mp^{i_t}} \right).
$$

\begin{lem} \label{l:nopquadro} 
The homology groups $H_*(\Br(*);\Z[q^\pmu])$ have no $p^2$--torsion.
\end{lem}

\begin{proof}
Notice that the monomials $1, x_1$ generate the groups $H_0(\Br(0), \Z[q^\pmu])$ and $H_0(\Br(1), \Z[q^\pmu])$ and that both these modules are equal to $\Z[q^\pmu]$.
For $i \geq 2$ the groups 
$H_0(\Br(i), \Z[q^\pmu]) = \Z$ are generated by the monomials $x_1^i$.
Now consider the following monomials:
\begin{gather} \label{eq:uno}
y_m^i, \quad y_m^ix_1 \qquad \mbox{ (case }p=0 \mbox{)}
\\ \label{eq:due}
\left.\begin{array}{l}
x_{2^j}^{2i}, \quad x_{2^j}^{2i}x_1, \\
x_{m2^j}^{2i}, \quad x_{m2^j}^{2i}x_1
\end{array}\right\}
\qquad  \mbox{ (case }p=2 \mbox{)}
\\ \label{eq:tre}
\left.\begin{array}{l}
y_{p^j}^{i}, \quad y_{p^j}^{i}x_1, \\
y_{2p^j}^{i}, \quad y_{2p^j}^{i}x_1, \\
y_{mp^j}^i, y_{mp^j}^ix_1
\end{array}\right\}
\qquad \mbox{ (case }p>2 \mbox{)}
\end{gather}
Because of the computations over $\Q[q^\pmu]$ (\cite{dps}, \cite{mar}), their liftings generate a free $\Z$--module of type $\Z[q^\pmu]/(\ph_h)$ in dimension $d$ in the homology of $\Br(n)$ whenever $n=kh$ or $n=kh+1$ and $d=k(h-2)$ and the Bockstein is zero for all these monomials.

All the other monomials lift to torsion classes and all these classes don't have $p^2$--torsion for any prime $p$. To prove this it is enough to show that in the submodule $M_p \sst H_*(\Br(*), \Z_p[q^\pmu])$ generated by all the monomials different from the ones in \eqref{eq:due} or \eqref{eq:tre}, we have that
$$
\ker \beta_p = \im \beta_p.
$$
Let us consider the set $S$ of the monomials that appear in the
polynomial rings of parts (a), (b) and (c) of
\fullref{p:phitorsione} and different from these in \eqref{eq:due} and
\eqref{eq:tre}.

Let us say that a monomial $w$ \emph{rises up} to a monomial $w'$ is $w$ appears as a summand in $\beta_p w'$. We call $w$ a \emph{basic} monomial if it doesn't appear as a summand in $\beta_p w'$ for any monomial $w'$. We also say that a monomial $w$ is a \emph{child} of $w'$ if $w'$ is basic and we can rise up from $w$ to $w'$ in a finite number of steps. We notice that in general a basic polynomial can be of the following form:
$$
w = x_{i_1}^{2a_1}\cdots x_{i_k}^{2a_k}y_{j_1}^{b_1}\cdots y_{j_h}^{b_k}x_{l_1} \cdots x_{l_s}.
$$
Let $\Delta_w$ be the set of all monomials that are children of $w$ (including $w$ itself). It is easy to see that $\Delta_w$ is in bijection with the set of the parts of $\{1, \ldots, s\}$ (with $s \geq 1$) if $l_1, \ldots, l_s$ are all different from $1$, or with the set of the parts of $\{1, \ldots, s-1 \}$ (with $s \geq 2$) if one of $l_1, \ldots, l_s$ is $1$.

Let we say that a monomial $w$ has $\ph$--torsion (over the ring $\Z_p[q^\pmu]$) if it generates a module isomorphic to $\Z_p[q^\pmu]/(\ph)$.
If a monomial $w$ has $\ph$--torsion over $\Z_p[q^\pmu]$ then all the other monomials children of $w$ have the same torsion. Moreover consider the algebraic complex $(M_w, \beta_p)$ given by the module $M_w$ generated (over $\Z_p[q^\pmu]$) by all the monomials in $\Delta_w$ and with the restriction of $\beta_p$ to $M_w$ as a boundary map: we have that $(M_w, \beta_p)$ is isomorphic to the algebraic complex that computes the reduced homology of the $(s-1)$--dimensional simplex with constant coefficients, over the ring $\Z_p[q^\pmu]/(\ph)$ and so $\ker \beta_p = \im \beta_p$ on $M_w$. 
\vspace{4pt}

One can check that for every monomial $w$ in $S$ there exists one and only one basic monomial $w'$ such that $w$ is a child of $w'$. This implies that the family of all different sets $\Delta_w$ gives a partition of $S$ and so $\ker \beta_p = \im \beta_p$ on all the module $M$. The Lemma follows. 
\end{proof}
\vspace{4pt}

\section{Main Result}
\vspace{4pt}

As a consequence of the last Lemma and of the previous computations, we can now state our main Theorem. Recall that the ring $ H_*(\Br(*); R[q^\pmu])$ is the bigraded direct sum of the modules $H_i(\Br(n); R[q^\pmu]) = H_i(\bd{F}_1(n),R)$, where $\bd{F}_1(n)$ is the Milnor fiber of the discriminant singularity for $\Br(n)$.
\vspace{4pt}

\begin{teo} \label{t:summary}

Set $\deg x_i= \deg y_i = i, \dim x_1 =0, \dim x_i = i-1$, $\dim y_i = i-2$. Then: 
$$ 
H_*(\Br(*); \Q[q^\pmu]) \simeq \Q[q^\pmu] \left[
\begin{array}{c}
x_1, y_m, m > 2
\end{array} 
\right] \left/ \left( 
\begin{array}{c}
\ph_2 x_1^2, \ph_m y_m
\end{array}
\right) \right. ;
$$
$$ 
H_*(\Br(*); \Z_2[q^\pmu]) \simeq \Z_2[q^\pmu] \left[
\begin{array}{c}
x_1, y_m, x_{2^{i+1}}^2, x_{m2^i}^2,\\
x_{2^i}^2 x_{2^{i_1}} \cdots x_{2^{i_h}}, \\
x_{m2^i}^2 x_{m2^{i_1}} \cdots x_{m2^{i_h}},  \\ 
y_m x_{m2^{j_1}} \cdots x_{m2^{j_h}}, \\
m \geq 2, 2 \nmid m
\end{array} 
\right] \left/ \left( 
\begin{array}{c}
\ph_2^{2^i} x_{2^i}^2, \\
\ph_m y_m, \\ 
\ph_m^{2^i} x_{m2^i}^2
\end{array}
\right) \right.
$$
with $0 \leq i$, $i+1 < i_1< \cdots < i_h$, 
     $0  < j_1 < \cdots < j_h$;

for $ p > 2 $: 
$$ 
H_*(\Br(*); \Z_p[q^\pmu]) \simeq \Z_p[q^\pmu] \left[
\begin{array}{c}
x_1, y_{p^i}, y_{mp^j}, y_{2p^i}\\
y_{p^i} x_{p^{i_1}} \cdots x_{p^{i_h}}, \\
x_1^2x_{2p^{j_1}} \cdots x_{2p^{j_h}} \\
y_{2p^i} x_{2p^{i_1}} \cdots x_{2p^{i_h}}, \\
y_{mp^j} x_{mp{j_1}} \cdots x_{mp^{j_h}} \\
m > 2, p \nmid m
\end{array} 
\right] \left/ \left( 
\begin{array}{c}
\ph_2 x_1^2, x_{2p^i}^2,\\
x_{p^i}^2, x_{mp^j}^2,\\
\ph_p^{\api} y_{p^i}, \\
\ph_2^{\bpi}y_{2p^i},\\ 
 \ph_m^{\phi(p^j)} y_{mp^j} 
\end{array}
\right) \right.
$$
with $0 < i < i_1< \cdots < i_h$, $0 \leq j < j_1 < \cdots < j_h$.
Finally, using the notation of the proof of \fullref{l:nopquadro}, we have:
$$
H_*(\Br(*);\Z[q^\pmu]) \simeq 
$$
$$
\Z[q^\pmu] \left[
\begin{array}{c}
x_1, y_m, m > 2
\end{array} 
\right] \left/ \left( 
\begin{array}{c}
\ph_2 x_1^2, \ph_m y_m
\end{array}
\right) \right. \oplus \bigoplus_{p\geq 2} \beta_p M_p. 
$$
\end{teo}

In tables \ref{ta:1}, \ref{ta:2}, \ref{ta:3} and \ref{ta:4} 
we give the explicit computations for some cases. The results in \fullref{ta:3} correspond to that in \cite{dps} for cohomology. We use the notation $\ph_h^i$ for the module $\Z_p[q]/(\ph_h^i)$ or $\Q[q]/(\ph_h^i)$ (note that $R[q^\pmu]/(\ph_h^i) = R[q]/(\ph_h^i)$). In \fullref{ta:4} we describe the additive structure of the integral homology of the fiber ${\bf F}_1(n)$.

\begin{table}[htb]
\begin{center}
\begin{tabular}{|c|c|c|c|c|c|c|c|c|c|}
\hline\strutt
$n$&$2$&$3$&$4$&$5$&$6$&$7$&$8$&$9$&$10$\\
\hline\strutt
$H_0$& $\ph_2$  &$\ph_2 $&$\ph_2 $&$\ph_2 $&$\ph_2 $&$\ph_2 $&$\ph_2 $&$\ph_2 $&$\ph_2 $\\
\hline\strutt
$H_1$ & & $\ph_3$ & $\ph_3$ & &&&&&\\
\hline\strutt
$H_2$ & & & $\ph_2^2$ & $\ph_2^2$ & $\ph_2 \oplus \ph_3$ & $\ph_2\oplus\ph_3$ & $\ph_2$ & $\ph_2$ & $\ph_2$\\
\hline\strutt
$H_3$ & & & & $\ph_5$ & $\ph_2 \oplus \ph_5$& $\ph_2$ &$\ph_2$ &$\ph_2 \oplus \ph_3$ &$\ph_2\oplus \ph_3$ \\
\hline\strutt
$H_4$ & & & & & $\ph_3$& $\ph_3$ & $\ph_2^2$ & $\ph_2^2$& $\ph_2$ \\
\hline\strutt
$H_5$ & & & & & & $\ph_7$ & $\ph_7$ & $\ph_3$ & $\ph_2 \oplus \ph_3$ \\
\hline\strutt
$H_6$ & & & & & & & $\ph_2^4$ & $\ph_2^4 \oplus \ph_3$ & $\ph_2 \oplus \ph_3 \oplus \ph_5$\\
\hline\strutt
$H_7$ & & & & & & & &$\ph_9$ & $\ph_2 \oplus \ph_9$\\
\hline\strutt
$H_8$ & & & & & & & & & $\ph_5$ \\
\hline
\end{tabular} 
\end{center}
\caption{$H_*(\Br(n); \Z_2[q^\pmu])$}
\label{ta:1}
\end{table}

\begin{table}[htb]
\begin{center}
\begin{tabular}{|c|c|c|c|c|c|c|c|c|c|}
\hline\strutt
$n$&$2$&$3$&$4$&$5$&$6$&$7$&$8$&$9$&$10$\\
\hline\strutt
$H_0$& $\ph_2$  &$\ph_2 $&$\ph_2 $&$\ph_2 $&$\ph_2 $&$\ph_2 $&$\ph_2 $&$\ph_2 $&$\ph_2 $\\
\hline\strutt
$H_1$ & & $\ph_3$ & $\ph_3$ & &&&&&\\
\hline\strutt
$H_2$ & & & $\ph_4$ & $\ph_4$ & $ \ph_3$ & $\ph_3$ &  &  & \\
\hline\strutt
$H_3$ & & & & $\ph_5$ & $ \ph_5$&  & & $\ph_3$& $\ph_3$\\
\hline\strutt
$H_4$ & & & & & $\ph_2^2$& $\ph_2^2$ & $\ph_2\oplus \ph_4$ & $\ph_2\oplus \ph_4$& $\ph_2$ \\
\hline\strutt
$H_5$ & & & & & & $\ph_7$ & $\ph_2 \oplus \ph_7$ & $\ph_2$ & $\ph_2$ \\
\hline\strutt
$H_6$ & & & & & & & $\ph_8$ & $\ph_8$ & $\ph_5$\\
\hline\strutt
$H_7$ & & & & & & & &$\ph_3^2$ & $\ph_3^2$\\
\hline\strutt
$H_8$ & & & & & & & & & $\ph_{10}$ \\
\hline
\end{tabular} 
\end{center}
\caption{$H_*(\Br (n); \Z_3[q^\pmu])$}
\label{ta:2}
\end{table}

\begin{table}[ht!]
\begin{center}
\begin{tabular}{|c|c|c|c|c|c|c|c|c|c|}
\hline\strutt
$n$&$2$&$3$&$4$&$5$&$6$&$7$&$8$&$9$&$10$\\
\hline\strutt
$H_0$& $\ph_2$  &$\ph_2 $&$\ph_2 $&$\ph_2 $&$\ph_2 $&$\ph_2 $&$\ph_2 $&$\ph_2 $&$\ph_2 $\\
\hline\strutt
$H_1$ & & $\ph_3$ & $\ph_3$ & &&&&&\\
\hline\strutt
$H_2$ & & & $\ph_4$ & $\ph_4$ & $\ph_3$ & $\ph_3$ &  &  & \\
\hline\strutt
$H_3$ & & & & $\ph_5$ & $\ph_5$&  & &$ \ph_3$ &$\ph_3$ \\
\hline\strutt
$H_4$ & & & & & $\ph_6$& $\ph_6$ & $ \ph_4$ & $ \ph_4$&  \\
\hline\strutt
$H_5$ & & & & & & $\ph_7$ & $\ph_7$ &  &  \\
\hline\strutt
$H_6$ & & & & & & & $\ph_8$ & $\ph_8 $ & $\ph_5$\\
\hline\strutt
$H_7$ & & & & & & & &$\ph_9$ & $\ph_9$\\
\hline\strutt
$H_8$ & & & & & & & & & $\ph_{10}$ \\
\hline
\end{tabular} 
\end{center}
\caption{$H_*(\Br (n); \Q[q^\pmu])$ }
\label{ta:3}
\end{table}

\begin{table}[htb]
\begin{center}
\begin{tabular}{|c|c|c|c|c|c|c|c|c|c|}
\hline\strutt
$n$&$2$&$3$&$4$&$5$&$6$&$7$&$8$&$9$&$10$\\
\hline\strutt
$H_0$& $\Z$  &$\Z $&$\Z $&$\Z $&$\Z $&$\Z $&$\Z $&$\Z $&$\Z $\\
\hline\strutt
$H_1$ & & $\Z^2$ & $\Z^2$ & &&&&&\\
\hline\strutt
$H_2$ & & & $\Z^2$ & $\Z^2$ & $\Z_2 \oplus \Z^2$ & $\Z_2\oplus \Z^2$ & $\Z_2$ & $\Z_2$ & $\Z_2$\\
\hline\strutt
$H_3$ & & & & $\Z^4$ & $\Z^4$&  & &$ \Z^2$ &$\Z^2$ \\
\hline\strutt
$H_4$ & & & & & $\Z^2$& $\Z^2$ & $\Z_{3}\oplus \Z^2$ & $\Z_3\oplus \Z^2$& $\Z_2 \oplus \Z_3$ \\
\hline\strutt
$H_5$ & & & & & & $\Z^6$ & $\Z^6$ & $\Z_2^2$ & $ \Z_2^2$ \\
\hline\strutt
$H_6$ & & & & & & & $\Z^4$ & $\Z^4 $ & $ \Z_2 \oplus \Z^4$\\
\hline\strutt
$H_7$ & & & & & & & &$\Z^6$ & $\Z^6$\\
\hline\strutt
$H_8$ & & & & & & & & & $\Z^4$ \\
\hline
\end{tabular} 
\end{center}
\caption{$H_*({\bd F}_1(n); \Z)$}
\label{ta:4}
\end{table}

Now consider the natural embeddings $j_n\co \Br(n) \into \Br(n+1)$. By definition the direct limit $\displaystyle \lim_{\longrightarrow}{}_n\Br(n) $ is the braid group on infinitely many strands $\Br(\infty)$. 
Notice that the first $p$--torsion class in the groups $H_*(\Br(n); \Z[q^\pmu])$ appears for $n = 2p + 2$ in dimension $2p-2$ and is stable; the corresponding generator is $x_1^2 x_2^2$ for $p = 2$ or $x_1^2 y_{2p}$ for $p > 2$. An equivalent result for the cohomology was proven in \cite{c}.

As a consequence of \fullref{t:summary}, we can compute the stable homology of braid groups, that is, the homology of $\Br(\infty)$ (see also \cite{cp}).

\begin{cor} \label{c:stable}
We have that
$
H_*(\Br(\infty); \Q[q^\pmu])= \Q,
$
concentrated in dimension $0$;
$$
H_*(\Br(\infty); \Z_2[q^\pmu])= \Z_2 [x_2^2, x_{2^i}, i>1]
$$
and for a prime $p > 2$
$$
H_*(\Br(\infty); \Z_p[q^\pmu])= \Z_p 
[y_{2p^i}, x_{2p^i}, i > 0 ]/(x_{2p^{i}}^2)
$$
with $\dim x_i = i-1$, $\dim y_i = i-2$. The Bockstein operator acts as follows:
$$
\beta_2 x_{2^i} = x_{2^{i-1}}^2;\qquad
\beta_p y_i = 0; \qquad \beta_p x_i = y_i \qquad \mbox{(for }p>2 \mbox{)}.
$$
The homology $H_*(\Br(\infty); \Z[q^\pmu])$ has no $p^2$--torsion for any prime $p$. A presentation of $H_*(\Br(\infty); \Z[q^\pmu])$ is given by
$$
\Z\left[
\begin{array}{c}
y_{2p^i},x_{2^j}^2,\\
x_{2^i}^2x_{2^{i_1}} \cdots x_{2^{i_h}} \\
y_{2p^j} x_{2p^{j_1}} \cdots x_{2p^{j_h}}
\end{array} 
\right] \left/ \left( 
\begin{array}{c}
2 x_{2^i}, py_{2p^j}, x_{2p^j}^2
\end{array}
\right) \right.
$$
with indices running as follows: $0<i$,$i+1 < i_1 < \cdots < i_h$, $0<j< j_1 < \cdots < j_h$ and $p$ in the set of odd primes. The structure of $\Z[q^\pmu]$--module is trivial and so the action of $q$ corresponds to multiplication by $-1$. 
\end{cor}

\appendix

\section*{Appendix} %\label{s:app}
Here we give the proof of a technical lemma stated in \fullref{s:milnor}.
\newtheorem*{lemn}{\fullref{l:tecnico1}}

\begin{lemn}
Let $m < n$ be two positive integers.
Also let $p$ be a prime. Then we have:
$$
(\ph_m, \ph_n) = \left\{ \begin{array}{cl}
(\ph_m, p) & \mbox{if }n=mp^i, i \geq1  \\ 
(1) & \mbox{otherwise.}
                         \end{array} \right.
$$
\end{lemn}

\begin{proof}
First of all, notice that the polynomials $\ph_m$ are irreducible for all $m \in \N$; hence, the quotient rings $\Z[q]/(\ph_m)$ are integral domains.

(i)\qua First suppose that $m \nmid n$ and let $l = \mbox{lcm}(m,n)$. Moreover we set $m' = \frac{l}{m}$, $n' = \frac{l}{n}$.
We have that $\ph_n \mid \frac{[l]}{[m]}$ and $\ph_m \mid \frac{[l]}{[n]}$.
Furthermore $\frac{[l]}{[m]} \cong m' \pmod{\ph_m}$ and $\frac{[l]}{[n]} \cong n' \pmod{\ph_n}$.
Since we have $(m',n') = (1)$ it follows that $(\ph_m,\ph_n) = (1)$. Hence the polynomial $\ph_n$ is invertible in $\Z[q]/(\ph_m)$ (and $\ph_m$ is invertible in $\Z[q]/(\ph_n)$).

(ii)\qua Now we suppose that $m \mid n$. For a fixed $m$ we want to prove by induction on $n$ that, modulo the multiplication by an invertible element in $\Z[q]/(\ph_m)$, the following holds:
\begin{eqnarray*}
\ph_n \cong p & \mbox{if }n=mp^i, i \geq 1\\
\ph_n \cong 1 & \mbox{otherwise.}
\end{eqnarray*}
If $n = mp$ we have that
$$
\frac{[n]}{[m]} \cong p \pmod{\ph_m}
$$
and so we can write
$$
\frac{[n]}{[m]} = \ph_n \prod _{
\begin{array}{l}
m' \mid m, m'<m\\ p \nmid m/m'
\end{array}
} \ph_{pm'}
$$
Since all the factors in the product are invertible, it follows that, modulo multiplication by invertible elements in $\Z[q]/(\ph_m)$, we get $\ph_n \cong p$.
\vspace{10pt}

If $n = mp^i$, in a similar way the next equality holds:
$$
\frac{[n]}{[m]}\ \ =\ \  \ph_n\!\!\!\!\!\! \prod _{\scriptsize	
\begin{array}{c}
1 \leq j <i \\ m' \mid m, m'p^j \neq n \\ p \nmid m/m'
\end{array}
}\!\!\!\!\!\! \ph_{m'p^j}\ \ \cong\ \ p^i \pmod{\ph_m}
$$
In the product there are exactly $i-1$ factors congruent to $p$ and all the others are invertible, so modulo invertible elements we have that $\ph_n \cong p$.
\vspace{10pt}

Finally we consider the case $n = mp_1^{i_1} \cdots p_k^{i_k}$. Let us define the set
$$
\cl{I} = \left\{
(m',j_1,\ldots,j_k) \in \N^{k+1}
\left|
\begin{array}{c}
m' \mid m,\\
0 \leq j_s \leq i_s \forall s,\\
(j_1,\ldots,j_k) \neq (0,\ldots,0),\\
m' \neq m \mbox{ if }j_s = i_s \forall s,\\
p_s \nmid (m/m') \quad \forall s \mbox{ s.t. } j_s \neq 0
\end{array}
\right\}.
\right.
$$
\vspace{10pt}
We have that:
%\begin{eqnarray*}
$$
\frac{[n]}{[m]} =  \prod _{n'\mid n, n' \nmid m} \ph_{n'} = 
   \prod_{I \in \cl{I}}\ph_{m'p_1^{j_1} \cdots p_k^{j_k}} \cong 
p_1^{i_1} \cdots p_k^{i_k} \pmod{\ph_m}
$$
%\end{eqnarray*}
and, by the inductive hypothesis, in the product, for all $s$ there are exactly $i_s$ factors congruent to $p_s$ and all the other factors are invertible; hence $\ph_n$ must be invertible, too. So the Lemma is proved. 
\end{proof}

\subsubsection*{Acknowledgments}
The author is grateful to M Salvetti for fruitful discussions. The
author is also thankful to F\,R Cohen and to the referee for many
useful suggestions.

\bibliographystyle{gtart}
\bibliography{link}

\end{document}